\input amstex
\documentstyle {amsppt}
\UseAMSsymbols \vsize 18cm \widestnumber\key{ZZZZZ}

\catcode`\@=11
\def\displaylinesno #1{\displ@y\halign{
\hbox to\displaywidth{$\@lign\hfil\displaystyle##\hfil$}&
\llap{$##$}\crcr#1\crcr}}
\def\ldisplaylinesno #1{\displ@y\halign{
\hbox to\displaywidth{$\@lign\hfil\displaystyle##\hfil$}&
\kern-\displaywidth\rlap{$##$} \tabskip\displaywidth\crcr#1\crcr}}
\catcode`\@=12

\refstyle{A}

\let \ol=\overline

\let \ti=\widetilde

\font\smain=cmsy10 at 7.5pt \font\ssmain=cmsy10 at 5.625pt

\font \fin=lasy8 at 15.4 pt

\def \so{\mathop{\hbox{\smain O}}\nolimits}
\def \sso{\mathop{\hbox{\ssmain O}}\nolimits}

\topmatter
\title  Standard Module Conjecture \endtitle

\rightheadtext{Standard Module Conjecture }
\author V. Heiermann and G. Mui\'c \endauthor

\address V.H.: Institut f\"ur Mathematik, Humboldt-Universit\"at zu
Berlin, Unter den Linden 6, 10099 Berlin, Allemagne\endaddress

\email heierman\@math.hu-berlin.de \endemail

\address G.M.: Department of Mathematics, University of Zagreb,
Bijenicka 30, 10000 Zagreb, Croatia\endaddress

\email gmuic\@math.hr \endemail

\abstract Let $G$ be a quasi-split $p$-adic group. Under the
assumption that the local coefficients $C_{\psi }$ defined with
respect to $\psi $-generic tempered representations of standard
Levi subgroups of $G$ are regular in the negative Weyl chamber, we
show that the standard module conjecture is true, which means that
the Langlands quotient of a standard module is generic if and only
if the standard module is irreducible.

\endabstract

\endtopmatter

\document

Let $F$ be a non archimedean local field of characteristic $0$.
Let $G$ be the group of points of a quasi-split connected
reductive $F$-group. Fix a $F$-Borel subgroup $B=TU$ of $G$ and a
maximal $F$-split torus $T_0$ in $T$. If $M$ is any semi-standard
$F$-Levi subgroup of $G$, a standard parabolic subgroup of $M$
will be a $F$-parabolic subgroup of $M$ which contains $B\cap M$.

Denote by $W$ the Weyl group of $G$ defined with respect to $T_0$
and by $w_0^G$ the longest element in $W$. After changing the
splitting in $U$, for any generic
representation $\pi $  of $G$, one can always find a non degenerate
character $\psi $ of $U$, which is compatible with $w_0^G$,
 such that $\pi $ is $\psi $-generic \cite{Sh2, section 3}. For any
semi-standard Levi-subgroup $M$ of $G$, we will still denote by
$\psi $ the restriction of $\psi $ to $M\cap U$. It is compatible
with $w_0^M$. If we write in the sequel that a representation of a
$F$-semi-standard Levi subgroup of $G$ is $\psi $-generic, then we
always mean that $\psi $ is a non degenerate character of $U$ with
the above properties.

Let $P=MU$ be a standard parabolic subgroup of $G$ and $T_M$ the
maximal split torus in the center of $M$. We will write $a_M^*$
for the dual of the real Lie-algebra of $T_M$ and $a_M^{*+}$ for
the positive Weyl chamber in $a_M^*$ defined with respect to $P$.
There is a canonical map $H_M:M\rightarrow a_M$, such that
$\vert\chi (m)\vert _F=q^{-\langle\chi ,H_M(m)\rangle}$ for every
$F$-rational character $\chi\in a_M^*$ of $M$. (We remark that this
is not the classical definition of $H_M$.) If $\pi $ is a smooth
representation of $M$ and $\nu \in a_M^*$, we denote $\pi _{\nu }$
the smooth representation of $M$ defined by $\pi _{\nu
}(m)=q^{-\langle\nu ,H_M(m)\rangle}\pi (m)$. The symbol $i_P^G$
will denote the functor of parabolic induction normalized such
that it sends unitary representations to unitary representations,
$G$ acting on its space by right translations.

Let $\tau $ be a generic irreducible tempered representation of
$M$ and $\nu\in a_M^{*+}$. Then the induced representation
$i_P^G\tau _{\nu }$ has a unique irreducible quotient $J(\tau ,\nu
)$, the so-called Langlands quotient.

The aim of our paper is to prove the \it standard module
conjecture \cite{CSh}, \rm which states that

\centerline{\it $J(\tau ,\nu )$ is generic, if and only if
$i_P^G\tau _{\nu }$ is irreducible. \rm}

\null We achieve this aim under the assumption that the local
coefficients $C_{\psi }$ defined with respect to $\psi $-generic
tempered representations of standard Levi subgroups of $G$ are
regular in the negative Weyl chamber. This property of the local
coefficients $C_{\psi }$ would be a consequence of Shahidi's
tempered $L$-function conjecture \cite{Sh2, 7.1}, which is now
known in most cases \cite{K}. Nevertheless, the result that we
actually need may be weaker (in particular, we do not need to
consider each component $r_i$ of the adjoint representation $r$
separately). So it may be possible to show it independently of the
tempered $L$-function conjecture (see the remark in {\bf 1.6}).

Our conditional proof of the standard module conjecture follows
the method developed in \cite{M}, \cite{M1}, but using the description of the
supercuspidal support of a discrete series representation of $G$
given in \cite{H1}.

The second named author would like to thank the Alexander von
Humboldt-foundation for their generous grant helping me to stay in
Berlin and complete this work.

\null\null {\bf 1.} Let $P=MU$ be a standard $F$-parabolic
subgroup of $G$ and $(\pi ,V)$ an irreducible $\psi $-generic
admissible representation of $M$. The parabolic subgroup of $G$
which is opposite to $P$ will be denoted $\ol{P}=M\ol{U}$. The set
of reduced roots of $T_M$ in $Lie(U)$ will be denoted $\Sigma
(P)$. We will use a superscript $G$ to underline that the
corresponding object is defined relative to $G$.

\null {\bf 1.1} For all $\nu $ in an open subset of $a_M^*$ we
have an intertwining operator $J_{\overline{P}\vert P}(\pi _{\nu
}):i_P^G\pi _{\nu }\rightarrow i_{\ol{P}}^G\pi _{\nu }$. For $\nu
$ in $(a_M^*)^+$ far away from the walls, it is defined by a
convergent integral $$(J_{\ol{P}\vert P}(\pi _{\nu })v)(g)=\int
_{\ol{U}} v(ug) du.$$ It is meromorphic in $\nu $ and the map
$J_{P\vert\ol{P}}J_{\ol{P}\vert P}$ is scalar. Its inverse equals
Harish-Chandra's $\mu $-function up to a constant and will be
denoted $\mu (\pi ,\nu )$.

\null {\bf 1.2} Put $\ti {w}=w_0^Gw_0^M$. Then
$\ti{w}\overline{P}\ti {w}^{-1}$ is a standard parabolic subgroup
of $G$. For any $\nu\in a_M^*$ there is a Whittaker functional
$\lambda _P(\nu, \pi ,\psi )$ on $i_P^GV$. It is a linear
functional on $i_P^GV$, which is holomorphic in $\nu $, such that
for all $v\in i_P^GV$ and all $u\in U$ one has $\lambda _P (\nu,
\pi ,\psi )((i_P^G\pi _{\nu })(u)v)=\ol{\psi (u)}\lambda
_{P}(\nu,\pi ,\psi )(v)$. Remark that by Rodier's theorem
\cite{R}, $i_P^G\pi _{\nu }$ has a unique $\psi $-generic
irreducible sub-quotient.

Fix a representative $w$ of $\ti {w}$ in $K$. Let $t(w)$ be the
map $i_{\overline{P}}^GV\rightarrow i_{w\ol{P}}^GwV$, which sends
$v$ to $v(w^{-1}\cdot )$. There is a complex number $C_{\psi }(\nu
,\pi ,w)$ such that $\lambda (\nu ,\pi ,\psi )=C_{\psi }(\nu ,\pi
,w)\lambda (w\nu ,w\pi ,\psi )t(w)J_{\ol{P}\vert P}(\pi _{\nu })$.
The function $a_M^*\rightarrow\Bbb C$, $\nu\mapsto C_{\psi }(\nu
,\pi ,w)$ is meromorphic.

The local coefficient $C_{\psi }$ satisfies the equality $C_{\psi
}(\cdot ,\pi ,w)C_{\psi }(w(\cdot ),w\pi ,w^{-1})=\mu (\pi, \nu )$
\cite{Sh1}.

\null{\bf 1.3} We will use the following criterion which follows
easily from the definitions and Rodier's theorem \cite{R}:

\null {\bf Proposition:} \it If $(\pi ,V)$ is an irreducible
tempered representation of $M$ and $\nu\in a_M^{*+}$, then the
Langlands quotient of the induced representation $i_P^G\pi_{\nu }$
is $\psi $-generic if and only if $\pi $ is $\psi $-generic and $C_{\psi }(\cdot
,\pi ,w)$ is regular in $\nu $. \rm

\null {\bf 1.4} For $\alpha\in\Sigma (P)$, put $w_{\alpha
}=w_0^{M_{\alpha }}w_0^M.$ With this notation, one has the
following version of the multiplicative formula for the local
coefficient $C_{\psi }(\cdot ,\pi ,\psi )$ \cite{Sh1, proposition
3.2.1}:

\null{\bf Proposition:} \it Let $Q=NV$ be a standard parabolic
subgroup of $G$, $N\subseteq M$, and $\tau $ an irreducible
generic representation of $N$, such that $\pi $ is a
sub-representation of $i_{Q\cap M}^M\tau $. Then one has
$$C^G_{\psi }(\cdot ,\pi ,w)=\prod _{\alpha\in\Sigma (Q)-\Sigma (Q\cap M)}
C^{N_{\alpha }}_{\psi }(\cdot ,\tau ,w_{\alpha }).$$

\null Proof: \rm It follows from \cite{Sh1, proposition 3.2.1}
that $$C^G_{\psi }(\cdot ,\pi ,w)=\prod _{\alpha\in\Sigma (P)}
C^{M_{\alpha }}_{\psi }(\cdot ,\tau ,w_{\alpha }).$$ Now fix
$\alpha\in\Sigma (P)$. As $J_{\ol{P}\cap M_{\alpha }\vert P\cap
M_{\alpha }}^{M_{\alpha }}(\pi _{\nu })$ (resp. $\lambda
^{M_{\alpha }}(\nu ,\pi ,\psi )$) equals the restriction of
$J_{\ol{Q}\cap M_{\alpha }\vert Q\cap M_{\alpha }}^{M_{\alpha
}}(\tau _{\nu })$ (resp. $\lambda ^{M_{\alpha }}(\nu ,\tau ,\psi
)$) to the space of $i_{P\cap M_{\alpha }}^{M_{\alpha }}\pi _{\nu
}\subseteq i_{Q\cap M_{\alpha }}^{M_{\alpha }}\tau _{\nu }$, it
follows that
$$C_{\psi }^{M_{\alpha }}(\nu ,\pi ,w_{\alpha })=C_{\psi }^{N_{\alpha }}(\nu ,\tau
,w_{\alpha }).$$ Applying the above product formula to the
expression on the right, one gets the required identity.
\hfill{\fin 2}

\null {\bf 1.5} Recall the tempered $L$-function conjecture
\cite{Sh, 7.1}: if $\sigma $ is an irreducible tempered
representation of $M$ and $M$ is a maximal Levi-subgroup of $G_1$,
then for every component $r_i$ of the adjoint representation $r$
the $L$-function $L(s,\sigma ,r_{G_1,i})$ is holomorphic for
$\Re(s)>0$.

\null {\bf Proposition:} \it Let $\pi $ be an irreducible generic
tempered representation of $M$. Assume that the tempered
L-function holds for $\pi $ relative to any $M_{\alpha }$,
$\alpha\in\Sigma (P)$.

Then $C_{\psi }(w(\cdot ),w\pi ,w^{-1})$ is regular in
$(a_M^*)^+$.

\null Proof: \rm Let $\lambda\in (a_M^*)^+$. Denote by $s_{\alpha
}\ti{\alpha }$, $s_{\alpha }\in\Bbb R$, the orthogonal projection
of $\lambda $ on $a_M^{M_{\alpha }*}$. Then $s_{\alpha }>0$. By
proposition {\bf 1.4} applied to $\tau =\pi $, $$C_{\psi
}(w(\lambda ),w\pi ,w^{-1})=\prod _{\alpha\in\Sigma (P)} C_{\psi
}^{M_{w\alpha }}(-s_{\alpha }\ti{w\alpha }, w\pi ,w_{w\alpha }).$$
By \cite{Sh, 3.11, 7.8.1 and 7.3},
$$C_{\psi }^{M_{w\alpha }}(s\ \ti{w\alpha }, w\pi ,w_{w\alpha
})=*\prod _i {L(1-is, w\pi ,r_{M_{w\alpha },i})\over L(is,w\pi
,r_{M_{w\alpha },i})},$$ where $*$ denotes a monomial in $q^{\pm
s}$.

Now, by assumption, $L(\cdot ,w\pi ,r_{M_{w\alpha },i})$ is
regular in $1+is_{\alpha }$. As $1/L(is_{\alpha },w\pi ,r_{M_{w\alpha },i})$
is polynomial, this proves the proposition. \hfill{\fin 2}

\null {\bf 1.6} \it Remark: \rm In fact, what is really needed to
prove the above proposition is a result that may be weaker than
the tempered $L$-function conjecture: suppose for simplicity that
$\pi $ is square integrable and choose a standard parabolic
subgroup $Q=NV$ of $G$ and a unitary supercuspidal representation
$\sigma $ of $N$, $N\subseteq M$, and $\nu\in a_N^*$, such that
$\pi $ is a sub-representation of $i_{Q\cap M}^M\sigma _{\nu }$.
Then, by {\bf 1.4},
$$C^G_{\psi }(s\ti{\alpha },\pi ,w)=\prod _{\alpha\in\Sigma (Q)-\Sigma (Q\cap M)}
C^{M_{\alpha }}_{\psi }(s\ti{\alpha }+\nu ,\sigma ,w_{\alpha }).$$
This is, up to a meromorphic function on the real axes, equal to
$$\prod _{\alpha\in\Sigma (Q)-\Sigma (Q\cap M)}{L(1-i_{\alpha
}s,\sigma _{\nu },r_{M_{w_{\alpha }},i_{\alpha }})\over L(i_{\alpha
}s,\sigma _{\nu },r_{M_{w_{\alpha }},i})}$$ where $i_{\alpha }\in\{1,2\}$
and $i_{\alpha }=i_{\beta }$ if $\alpha $ and $\beta $ are
conjugated.

Now let $\Sigma _{\so }(N)$ be the set of reduced roots $\alpha $
of $T_N$ in $Lie(V)$, such that Harish-Chandra's $\mu $-function
$\mu ^{N_{\alpha }}$ defined with respect to $N_{\alpha }$ and
$\sigma $ has a pole. The set of these roots forms a root system
\cite{Si, 3.5}. Denote by $\Sigma _{\so }(Q)$ the subset of those
roots, which are positive for $Q$ and by $\Sigma _{\so }(Q\cap M)$
the one of those roots in $\Sigma _{\so }(Q)$, which belong to
$M$. Then the above product equals up to a holomorphic function
$$\prod _{\alpha\in\Sigma _{\sso}(Q)-\Sigma _{\sso }(Q\cap M)}{1-q^{-i_{\alpha
}(s+\langle\alpha^{\vee },\nu \rangle)}\over 1-q^{-1+i_{\alpha
}(s+\langle\alpha^{\vee },\nu \rangle)}},$$ So, what we have to
know, is that this meromorphic function is holomorphic for $s<0$.
The fact that $\sigma _{\nu }$ lies in the supercuspidal support
of a discrete series means by the main result of \cite{H1}, that
$\sigma _{\nu }$ is a pole of order $rk_{ss}(M)-rk_{ss}(N)$ of
Harish-Chandra's $\mu $-function $\mu ^M$. This can be translated
somehow to the assertion that $\nu $ corresponds to a
distinguished unipotent orbit \cite{H2}. So, what remains, is a
purely combinatorial problem in the theory of rootal hyperplane
configurations, which can be stated independently of
representation theory, although its validity may depend on the
fact that the labels $i_{\alpha }$ and $\nu $ come from a "generic
setting". As, in particular, one does not need here to consider
each component $r_i$ of the adjoint representation $r$ separately,
this result may be weaker than Shahidi's tempered $L$-function
conjecture.

\null \null {\bf 2.} In this section we make the following
assumption on $G$ (see the remark in {\bf 1.5} for what we
actually need):

\null \it (TL) If $M$ is a semi-standard Levi subgroup of $G$ and
if $\pi $ is an irreducible generic tempered representation of $M$
then $L(s,\pi ,r_i)$ is regular for $\Re(s)>0$ for every $i$. \rm

\null We give a proof of the following lemma only for
completeness:

\null {\bf 2.1. Lemma:} \it Let $P=MU$ be a $F$-standard parabolic
subgroup of $G$ and $\sigma $ an irreducible supercuspidal
representation of $M$. If the induced representation $i_P^G\sigma
$ has a sub-quotient, which lies in the discrete series of $G$,
then any tempered sub-quotient of $i_P^G\sigma $ lies in the
discrete series of $G$.

\null Proof: \rm If $i_P^G\sigma $ has an irreducible
sub-quotient, which is square-integrable, then by the main result
of \cite{H}, $\sigma $ is a pole of order $rk_{ss}(G)-rk_{ss}(M)$
of $\mu $ and $\sigma _{\vert A_M}$ is unitary. It follows that
the central character of $i_P^G\sigma $ is unitary, too, which
implies that the central character of any irreducible sub-quotient
of $i_P^G\sigma $ is unitary. In particular, any essentially
tempered irreducible sub-quotient of $i_P^G\sigma $ is tempered.

So, if $\tau $ is an irreducible tempered sub-quotient of
$i_P^G\sigma $, then there is a $F$-parabolic subgroup $P'=M'U'$
of $G$ and a square-integrable representation $\sigma '$ of $M'$,
such that $\tau $ is a sub-representation of $i_{P'}^G\sigma '$.
The supercuspidal support of $\sigma '$ and the $W$-orbit of
$\sigma $ share a common element $\sigma _0$.

By the invariance of Harish-Chandra's $\mu $-function, $\sigma _0$
is still a pole of $\mu $ equal to the order given above. As this
order is maximal and the central character of $\sigma '$ must be
unitary, this implies that $M'$ must be equal to $G$, and
consequently $\tau =\sigma '$ is square-integrable. \hfill{\fin 2}

\null {\bf 2.2 Theorem:} \it Let $G$ be a group that satisfies
property (TL). Let $P=MU$ be a $F$-standard parabolic subgroup of
$G$ and $\sigma $ be an irreducible $\psi $-generic supercuspidal
representation of $M$.

If the induced representation $i_P^G\sigma$ has a sub-quotient,
which lies in the discrete series of $G$ (resp. is tempered), then
any irreducible $\psi $-generic sub-quotient of $i_P^G\sigma $
lies in the discrete series of $G$ (resp. is tempered).

\null Proof: \rm First assume that $i_P^G\sigma $ has a
sub-quotient, which lies in the discrete series of $G$. Let $(\pi
,V)$ be an irreducible, admissible $\psi $-generic representation of $G$,
which is a sub-quotient of $i_P^G\sigma$. By the Langlands
quotient theorem, there is a standard parabolic subgroup
$P_1=M_1U_1$ of $G$, an irreducible tempered representation $\tau
$ of $M_1$ and $\nu \in (a_{M_1}^*)^+$, such that $\pi $ is the
unique irreducible quotient of $i_{P_1}^G\tau _{\nu }$.

As any representation in the supercuspidal support of $\tau _{\nu
}$ must lie in the supercuspidal support of $\pi $, any such
representation must be conjugated to $\sigma$. So, after
conjugation by an element of $G$, we can assume that $M\subseteq
M_1$ and that $\tau _{\nu }$ is a sub-representation of $i_{P\cap
M_1}^{M_1}\sigma$.

We will actually show that $P_1=G$, which means that $\pi $ is
tempered and by {\bf 2.1} in fact square-integrable.

Following {\bf 1.3}, $\tau $ must be $\psi $-generic and it is enough to
show that $C_{\psi }(\cdot ,\tau,w)$ has a pole in $\nu $, if
$P_1\ne G$.

For this we will use the assumption that $i_P^G\sigma $ has an
irreducible sub-quotient which is square-integrable. By the main
result of \cite{H} this implies that $\mu $ has a pole of order
equal to $rk_{ss}G-rk_{ss}M$ in $\sigma$. Remark that $\mu ^{M_1}$
can have at most a pole of order $rk_{ss}M_1-rk_{ss}M$ in
$\sigma$. The order of the pole of $\mu $ in $\tau _{\nu }$ is
equal to the one of $\mu /\mu ^{M_1}$ in $\sigma$. It follows that
the order of this pole must be $>0$, if $P_1\ne G$. As
$$\mu (\tau ,\cdot)=C_{\psi }(\cdot ,\tau ,w)C_{\psi }(w(\cdot ),w\tau, w^{-1}),$$
it follows that either $C_{\psi }(\cdot ,\tau ,w)$ or $C_{\psi
}(w(\cdot ),w\tau, w^{-1})$ must have a pole in $\nu $. As $\nu
\in (a_M^*)^+$, it follows from {\bf 1.5} that $C_{\psi }(w(\cdot
),w\tau, w^{-1})$ cannot have a pole in $\nu $. So $C_{\psi
}(\cdot ,\tau ,w)$ does. This gives us the desired contradiction.

Now assume that $i_P^G\sigma $ only has a tempered sub-quotient
$\tau $. Then there is a standard parabolic subgroup $P_1=M_1U_1$
of $G$ and a discrete series representation $\pi _1$ of $M_1$,
such that $\tau $ is a sub-representation of $i_{P_1}^G\pi _1$. As
the supercuspidal support of $\pi _1$ is contained in the
$G$-conjugacy class of $\sigma $, it follows that there is a
standard Levi subgroup $M'\supseteq M$, such that $i_{P\cap
M'}^{M'}\sigma $ has a discrete series sub-quotient.

By what we have just shown, there exists a unique $\psi $-generic
subquotient $\pi '$ of $i_{P\cap M'}^{M'}\sigma $, which lies in
the discrete series.

As $i_P^G\sigma $ and $i_{P'}^G\pi'$ have each one a unique
irreducible $\psi $-generic sub-quotient and any sub-quotient of
$i_{P'}^G\pi'$ is a sub-quotient of $i_P^G\sigma $, these
irreducible $\psi $-generic sub-quotients must be equal and
therefore tempered. \hfill{\fin 2}

\null {\bf 2.3 Theorem:} \it Let $G$ be a group that satisfies
property (TL). Let $P=MU$ be a $F$-standard Levi subgroup of $G$,
$\tau $ an irreducible tempered generic representation of $M$ and
$\nu\in a_M^{*+}$.

Then the Langlands quotient $J(\tau ,\nu )$ is generic, if and
only if $i_P^G\tau _{\nu }$ is irreducible.

\null Proof: \rm As $i_P^G\tau _{\nu }$ always has a generic
sub-quotient, one direction is trivial. So, assume $i_P^G\tau
_{\nu }$ is reducible. We will show that $\pi =J(\nu ,\tau )$ is
not $\psi $-generic for any $\psi $.

We can consider (and will) $\nu_{\pi}:=\nu$ as an element of
$a_T^{*}$. We denote by $<$ the partial order on $a_T^{*}$
explained in ([BW], Chapter XI, 2.1) (for our purpose it is not
important to write it explicitly).

Let $\pi '$ be an irreducible sub-quotient of $i_P^G\tau _{\nu }$,
which is not isomorphic to $\pi $. Let $P'=M'U'$ be a $F$-standard
parabolic subgroup, $\tau '$ an irreducible tempered
representation of $M'$ and $\nu '\in a_{M'}^{*+}$, such that $\pi
'=J(\nu ',\tau ')$. Let $\nu_{\pi'}:=\nu'$. Then \cite{BW, XI,
Lemma 2.13}
$$
\nu_{\pi'}<\nu_\pi. \tag 2.1
$$

Choose an $F$-standard parabolic
subgroup $P_1=M_1U_1$, $M_1\subseteq M$, with an irreducible
$\psi$--generic supercuspidal representation $\sigma $ of $M_1$, such that
$\tau $ is a sub-quotient of $i_{P_1\cap M}^M\sigma $.

Then $\sigma _{\nu }$ lies as well in the supercuspidal support of
$\pi $ as in the supercuspidal support of $\pi '$. It lies also in
the $G$-conjugacy class of the supercuspidal support of $\tau
'_{\nu '}$ and $\tau _{\nu }$. Let $\pi _{\psi }$ be the unique
$\psi $-generic irreducible sub-quotient of $i_{P_1}^G\sigma _{\nu
}$. By {\bf 2.2}, the unique $\psi $-generic irreducible
sub-quotient $\tau ''$ of $i_{P_1\cap M'}^{M'}\sigma $ is
tempered. The induced representation $i_{P'}^G\tau ''_{\nu '}$
admits a unique $\psi $-generic irreducible sub-quotient, which is
equal to the unique $\psi $-generic sub-quotient of
$i_{P_1}^G\sigma _{\nu }$.   Let
$\pi{''}=J(\nu', \tau '')$ be the Langlands quotient of $i_{P'}^G\tau ''_{\nu '}$. Since (2.1) implies
$\nu_{\pi{''}}=\nu_{\pi'}<\nu_\pi$, $\pi $ cannot be a
sub-quotient of $i_{P'}^G\tau ''_{\nu '}$ by \cite{BW, XI, Lemma 2.13}.
Therefore, $\pi $ is not $\psi $-generic.

\hfill{\fin 2}

\enddocument
\bye